\documentclass[a4paper,12pt]{article}

\newcommand{\qed}{\hfill $\sqcap\!\!\!\!\sqcup$} 

\newcommand{\noa}{\noalign{\medskip}}

\newcommand{\sgn}{\hbox{sgn}}

\newcommand{\pa}{\partial}

\newcommand{\di}{\displaystyle}

\begin{document}
\pagestyle{plain}
\title{Riemannian optimal control}

\author{Andreea Bejenaru and Constantin Udri\c{s}te}
\date{}
\maketitle

\begin{abstract} The aim of this paper is to adapt the general multitime maximum
principle to a Riemannian setting. More precisely, we intend to study geometric 
optimal control problems constrained by the metric compatibility evolution PDE system; the
evolution ("multitime") variables are the local coordinates on a
Riemannian manifold, the state variable is a Riemannian structure
and the control is a linear connection compatible to the Riemannian
metric. We apply the obtained results in order to
solve two flow-type optimal control problems on Riemannian setting: firstly,
we maximize the total divergence of a fixed vector field;
secondly, we optimize the total Laplacian (the gradient flux) of a fixed differentiable function.
Each time, the result is a bang-bang-type optimal linear connection.
Moreover, we emphasize the possibility of choosing at least two
soliton-type optimal (semi-) Riemannian structures. Finally, these theoretical examples 
help us to conclude about the geometric optimal shape of pipes, 
induced by the direction of the flow passing through them.

{\bf Keywords:} multitime maximum principle, Riemannian optimal control, shape optimization, 
gradient flow, total divergence, total Laplacian, bang-bang-type optimal solution, soliton-type metric.

{\bf MSC2010:}  49J20,  49N05, 49Q10, 53C05, 53C80.
\end{abstract}

\section{Adjoint PDE systems \\in Riemannian geometry}

A connection on a manifold is a type of differentiation that acts on
vector fields, differential forms and tensor products of these
objects (see \cite{[2]}, \cite{[6]}, \cite{[11]}). Its importance lies in the fact that given a piecewise
continuous curve connecting two points on the manifold, the
connection defines a linear isomorphism between the respective
tangent spaces at these points. Another fundamental concept in the
study of differential geometry is that of a Riemannian metric. It is
well known that a Riemannian metric uniquely determines a
Levi-Civita connection: a symmetric connection for which the
Riemannian metric is parallel. Since we may define linear connections unattached to 
Riemannian metrics, it is natural to ask, for a symmetric
connection, if there exists a parallel Riemannian metric, that is,
whether the connection is a Levi-Civita one. More generally, a
connection on a manifold $M$, symmetric or not, is said to be metric
if admits a parallel Riemannian metric defined on $M$. Then, the equations 
describing the metric property of symmetric linear connections are called 
{\it metric compatibility evolution PDE system}. 

Let $M$ be an $n$-dimensional differentiable manifold with local coordinates $(x^1,...,x^n)$. 
As we have mentioned above, fixing a Riemannian structure $g$ on $M$
ensures us about the existence of a symmetric linear connection satisfying the 
{\it metric compatibility PDE system} $\nabla_{\frac{\pa}{\pa x^i}}g=0,\,\,\forall i=1,...,n$.

Let us change the geometric point of view in those of deformation theory.
In this sense, let us discuss about the Riemannian metric controlled by a connection. 
For that we consider the {\it controlled evolution law} (linear PDE system)
$$\frac{\pa g_{ij}}{\pa x^k}(x)=g_{ps}(x)\left[\delta^p_i\Gamma^s_{jk}(x)+\delta^p_j\Gamma^s_{ik}(x)\right],\,\,i,j,k=1,...,n,\leqno(PDE)$$
together with the initial condition
$$g_{ij}(x_0)=\eta_{ij},\leqno(x_0)$$
where the piecewise metric tensor $g=(g_{ij})$ denotes a {\it symmetric
state tensor}, $x=(x^1,...,x^n)$ is the {\it multitime variable} (see \cite{[1]}, \cite{[3]},\cite{[12]}-\cite{[21]}),
and $\Gamma=(\Gamma^k_{ij})$ denotes the {\it symmetric control
linear connection}.

The PDE system has solutions if and only if the complete integrability conditions
$$\frac{\pa }{\pa x^l}\left\{g_{ps}\left[\delta^p_i\Gamma^s_{jk}+\delta^p_j\Gamma^s_{ik}\right]\right\}
=\frac{\pa }{\pa
x^k}\left\{g_{ps}\left[\delta^p_i\Gamma^s_{jl}+\delta^p_j\Gamma^s_{il}\right]\right\},\,\,\forall
i,j,k,l=1,...,n\leqno(CIC)$$ are satisfied. Explicitly, this
means $R_{ijkl}+R_{jikl}  = 0,$ where $R_{ijkl}$ denotes the Riemann curvature
tensor field corresponding to the solution $(g,\Gamma)$.

We consider the set of {\it admissible controls}
$$\mathcal{U}=\{\Gamma:M\rightarrow R^{n^3}|\,\,\Gamma_{ij}^k=\Gamma_{ji}^k\}.$$

Since the PDE system is linear, it coincides with its {\it infinitesimal deformation},
around a solution $g_{ij}(x)$. This PDE is also {\it auto-adjoint}
since $vp_{t^2} - p v_{t^2} = 0$, for any two solutions $v(x,t)$ and $p(x,t)$. 
If it is used as adjoint equation, then a solution $p(x, t)$ is called the {\it costate function}. 

The foregoing PDE systems determine the {\it multitime adjoint PDEs}
$$\frac{\pa \lambda^{ijk}}{\pa x^k}(x)=-\lambda^{psk}(x)[\delta^i_p\Gamma^j_{sk}(x)+\delta^i_s\Gamma^j_{pk}(x)],\leqno(ADJ)$$
whose solution $\lambda = (\lambda^{ijk})$, called the {\it costate
tensor}, is not necessary symmetric. The systems (PDE) and (ADJ) are
adjoint (dual) in the sense of zero divergence of the vector field
$Q(x)=(Q^k(x)=y_{ij}(x)\lambda^{ijk}(x)),$ where $y(x)=(y_{ij}(x))$ denotes an infinitesimal deformation
around a solution $g_{ij}(x)$.

The symmetry of the state and control variables suggests us to
consider a {\it symmetrized adjoint PDE system}, corresponding to
the {\it symmetric costate variables}
$p^{ijk}=\lambda^{ijk}+\lambda^{jik}.$ By computation, $p$ is
solution for
$$\frac{\pa p^{ijk}}{\pa x^k}(x)=-p^{rsk}(x)[\delta^i_r\Gamma^j_{sk}(x)+\delta^j_r\Gamma^i_{sk}(x)].\leqno(ADJ^s)$$
Again, $(PDE)$ and $(ADJ^s)$ are dual in the sense of zero
divergence of the vector field $S(x)=(S^k(x)=y_{ij}(x)p^{ijk}(x))$.

Moreover, introducing the Hamiltonian
$H(x,g,\Gamma,p)=p^{ijk}g_{is}\Gamma^s_{jk}$, the evolution systems
$(PDE)$ and $(ADJ^s)$ become
$$\frac{\pa g_{ij}}{\pa x^{k}}(x) =\left[\frac{\pa H}{\pa p^{ijk}}+\frac{\pa H}{\pa p^{jik}}\right](x,g(x),\Gamma(x),p(x)),\leqno(PDE)$$
respectively
$$\frac{\pa p^{ijk}}{\pa x^{k}}(x) = -\left[\frac{\pa H}{\pa g_{ij}}+\frac{\pa H}{\pa g_{ji}}\right](x,g(x),\Gamma(x),p(x)).\leqno(ADJ^s)$$


\section{Geometric optimal control with \\metric evolution-type constraints}

\subsection{Multitime maximum principle}

In this section we consider a general control problem, with
functional defined as multiple integral and evolution described by
the Riemannian metric compatibility PDE system.  Its solution is based on a general multitime maximum principle 
analyzed, under different aspects in \cite{[3]} and \cite{[12]}-\cite{[21]}, and generalizing 
the classical approach on single-time Pontriaguine maximum principle (see \cite{[4]},\cite{[7]}-\cite{[10]}, \cite{[22]}). More precisely,
in \cite{[21]} we have proved, using needle-shaped control variations, that, given the hyper-parallelepiped $\Omega_{0t_{0}}$ in $R^m$ having $0$ and $t_{0}$ as diagonal points, the solutions for the multitime
optimal control problem with running cost and integral cost on
boundary,
$$
\begin{array}{cc}
\hspace{0.cm}\di{\max_{u(\cdot)}}&\hspace{-1.3cm}\left(J[u(\cdot)]=\di\int_{\Omega_{0t_{0}}}X(t,x(t),u(t))\,dt+\di\int_{\partial\Omega_{0t_0}}\chi(t,x(t))\,
d\sigma\right) \\ \noa \hbox{subject to} \,\,&\di\frac{\pa x^i}{\pa
t^\alpha}(t) = X^i_\alpha (t,x(t),u(t)),\,\, i=1,...,n,\,\, \alpha =
1,..., m, \\ \noa &\hspace{-1.3cm}u(t)\in R^k, \,\,x(t)\in
R^n,\,\,t\in {\Omega_{0t_0}}\subset R^m, \,\,x(0)=x_0,
\end{array}
$$
 with the corresponding control Hamiltonian
$$H(t,x,p,u)=X(t,x,u)+p^{\alpha}_{i}X^{i}_{\alpha}(t,x,u),$$
satisfy the multitime maximum principle described in the next
fundamental outcome.

{\bf Theorem (Multitime maximum principle).} {\it Suppose
$u^*(\cdot)$ is an optimal solution of the control problem and
 $x^*(\cdot)$ is the corresponding optimal state. Then there exists
a costate tensor $p^*=(p^{*\alpha}_{i}): \Omega_{0t_0} \to R^{mn}$
such that

$$
\frac{\partial x^{*i}}{\partial t^{\alpha}}(t) =\frac{\partial
H}{\partial p^{\alpha}_{i}}(t,x^*(t),p^*(t),u^*(t)) ,
$$
$$
\frac{\partial p^{*\alpha}_{i}}{\partial t^{\alpha}}(t) =
-\frac{\partial H}{\partial x^{i}}(t,x^*(t),p^*(t),u^*(t)),
$$
$$H(t,x^*(t),p^*(t),u^*(t))=\di{\max_{u(\cdot)}}\,H(t,x^*(t),p^*(t),u(t)),\,\,\forall
t\in \Omega_{0t_{0}}$$ and
$$
n_{\alpha}p^{\ast\alpha}_{i}|_{\partial\Omega_{0t_{0}}}=\frac{\partial
\chi}{\partial x^{i}}|_{\partial\Omega_{0t_{0}}}.
$$}

\subsection{Riemannian optimal control}

If $(M,g)$ is a Riemannian manifold, let $x=(x^1,...,x^n)$ denote
the local coordinates relative to a fixed local map $(V,h)$. Since
$h:V\rightarrow R^n$ is an isomorphism, we denote by
$\Omega_{x_{0}x_{1}}$ a subset of $V$ diffeomorphic through $h$ with
the hyper-parallelepiped in $R^n$ having $h(x_{0})$ and $h(x_1)$ as
diagonal points.   If $X=X(x,g,\Gamma)$ and $\chi=\chi(x,g)$ are
differentiable maps, we associate the Bolza-type cost functional
$$J[\Gamma]=\int_{\Omega_{x_0x_1}}X(x,g(x),\Gamma(x))dx+\int_{\pa\Omega_{x_0x_1}}\chi(x,g(x))d\sigma,$$
where $dx$ denotes the differential $n$-form $dx^1\wedge...\wedge dx^n$ (the Euclidean volume element on 
$\Omega_{x_0x_1}$) and $d\sigma$ is the Euclidean volume element on the boundary.

The multitime optimal control problem consists in finding

$$\di{\max_\Gamma}\,\left(J[\Gamma]=\int_{\Omega_{x_0x_1}}X(x,g(x),\Gamma(x))\,dx+\int_{\pa\Omega_{x_0x_1}}\chi(x,g(x))\,d\sigma\right),\leqno(J)$$
subject to the evolution system
$$\frac{\pa g_{ij}}{\pa x^k}(x)=g_{ps}(x)\left[\delta^p_i\Gamma^s_{jk}(x)+\delta^p_j\Gamma^s_{ik}(x)\right],\,\,i,j,k=1,...,n\leqno(PDE)$$
and the initial condition
$$g(x_0)=\eta.\leqno(x_0)$$

Since the main ingredients of this Riemannian optimal control
problem (the state variables, the control variables and evolution
constraints $(PDE)$) are symmetric, we shall derive an adapted
multitime maximum principle, based on symmetric costate variables.
For this, we introduce the symmetric Lagrange multipliers
$p^{ijk}=p^{jik}$ and the {\it reduced control Hamiltonian}
$$H(x,g,\Gamma,p)=X(x,g,\Gamma)+g_{is}\Gamma^s_{jk}p^{ijk}.\leqno(H)$$

{\bf Corollary 1 (Riemannian maximum principle).} {\it Suppose the symmetric connection 
$\Gamma^*(\cdot)$ is an optimal solution for $((PDE),(J),(x_0))$ and that $g^*(\cdot)$ is
the corresponding optimal Riemannian structure. Then there exists
a symmetric dual tensor $p^*=(p^{*ijk}): \Omega_{x_0x_1} \to R^{n^3}$ such that
$$\frac{\partial g^*_{ij}}{\partial x^{k}}(x) =\left[\frac{\partial H}{\pa p^{ijk}}
+\frac{\partial H}{\pa p^{jik}}\right](x,g^*(x),\Gamma^\ast(x),p^*(x)),\leqno(PDE)$$

$$\frac{\pa p^{\ast ijk}}{\pa x^{k}}(x) = -\left[\frac{\pa H}{\pa g_{ij}}
+\frac{\pa H}{\pa g_{ji}}\right](x,g^*(x),\Gamma^\ast(x),p^*(x))\leqno(ADJ^s)$$
and
$$H(x,g^*(x),\Gamma^\ast(x),p^*(x))=\di{\max_{\Gamma(\cdot)}}\,H(x,g^*(x),\Gamma(x),p^*(x)),\,\,\forall x\in \Omega_{x_0x_{1}}.\leqno(Opt)$$
Finally, the boundary conditions
$$n_{k}p^{\ast ijk}|_{\partial\Omega_{x_0x_{1}}}=\left[\frac{\partial \chi}{\partial g_{ij}}
+\frac{\partial \chi}{\partial g_{ji}}\right]_{\partial\Omega_{x_0x_{1}}} \leqno (\partial\Omega_{x_0x_{1}})$$
are satisfied, where $n$ denotes the covector corresponding to the unit normal vector on $\partial\Omega_{x_0x_{1}}$.}

{\bf Proof.} We denote by $\overline{H}$ the standard control Hamiltonian 
corresponding to the Riemannian optimal control problem $((J),(PDE),(x_0))$, that is
$$\overline H(x,g,\Gamma,\lambda)=X(x,g,\Gamma)+\lambda^{ijk}\left[g_{is}\Gamma^s_{jk}(x,g,\Gamma)+g_{js}\Gamma^s_{ik}(x,g,\Gamma)\right]$$
$$=X(x,g,\Gamma)+g_{is}\Gamma^s_{jk}\left[\lambda^{ijk}+\lambda^{jik}\right]=H(x,g,\Gamma,\lambda^{ijk}+\lambda^{jik}).$$

Let us define the symmetric costate tensor
$p^{ijk}=\lambda^{ijk}+\lambda^{jik}$. Writing the multitime maximum
principle with standard Hamiltonian, and using the definition of
$p$, we obtain
$$\frac{\partial g^*_{ij}}{\partial x^{k}} =\frac{\partial \overline H}{\pa \lambda^{ijk}}=\frac{\partial H}{\pa p^{ijk}}+\frac{\partial H}{\pa p^{jik}};$$
$$\frac{\pa p^{\ast ijk}}{\pa x^{k}} =\frac{\pa \lambda^{\ast ijk}}{\pa x^{k}}+\frac{\pa \lambda^{\ast jik}}{\pa x^{k}} = -\left[\frac{\pa \overline H}{\pa g_{ij}}+\frac{\pa \overline H}{\pa g_{ji}}\right]=
-\left[\frac{\pa H}{\pa g_{ij}}+\frac{\pa H}{\pa g_{ji}}\right];$$
$$H(x,g^*,\Gamma^\ast,p^*)=\overline H(x,g^*,\Gamma^\ast,\lambda^*)=\di{\max_{\Gamma}}\,\overline H(x,g^*,\Gamma,\lambda^*)=\di{\max_{\Gamma}}\,H(x,g^*,\Gamma,p^*);$$
$$n_{k}p^{\ast ijk}|_{\partial\Omega_{x_0x_{1}}} =n_{k}\left[\lambda^{\ast ijk}+\lambda^{\ast jik}\right]_{\partial\Omega_{x_0x_{1}}}
=\left[\frac{\partial \chi}{\partial g_{ij}}+\frac{\partial
\chi}{\partial g_{ji}}\right]_{\partial\Omega_{x_0x_{1}}}.$$ \qed

{\bf Remark.} By replacing the metric compatibility evolution $(PDE)$ with the PDE system corresponding to the dual tensor $g^{-1}$:
$$\frac{\pa g^{ij}}{\pa x^k}(x)=-g^{ps}(x)\left[\delta^i_p\Gamma^j_{sk}(x)+\delta^j_p\Gamma^i_{sk}(x)\right],\,\,i,j,k=1,...,n,\leqno(PDE')$$
with initial condition
$$g^{ij}(x_0)=\eta^{ij} \leqno(x'_0)$$
and using the dual Hamiltonian 
$$H'(x,g^{-1},\Gamma,p)=X(x,g,\Gamma)-g^{is}\Gamma^j_{sk}p^{k}_{ij},\leqno(H')$$
we can rephrase the Riemannian multitime maximum principle as it follows.

{\bf Corollary 2 (Riemannian dual maximum principle).} {\it Suppose
the symmetric connection $\Gamma^*(\cdot)$ is an optimal solution
for $((PDE'),(J),(x'_0))$ and that $g^{*-1}(\cdot)$ is the
corresponding optimal state. Then there exists a symmetric dual
tensor $p^*=(p^{*k}_{ij}): \Omega_{x_0x_1} \to R^{n^3}$ such that
$$\frac{\partial g^{*ij}}{\partial x^{k}}(x) =\left[\frac{\partial H'}{\pa p^{k}_{ij}}+\frac{\partial H'}{\pa p^{k}_{ji}}\right](x,g^{*-1}(x),\Gamma^\ast(x),p^*(x)),\leqno(PDE')$$
$$\frac{\pa p^{\ast k}_{ij}}{\pa x^{k}}(x) = -\left[\frac{\pa H'}{\pa g^{ij}}+\frac{\pa H'}{\pa g^{ji}}\right](x,g^{*-1}(x),\Gamma^\ast(x),p^*(x))\leqno(ADJ'^s)$$
and
$$H'(x,g^{*-1}(x),\Gamma^\ast(x),p^*(x))=\di{\max_{\Gamma(\cdot)}}\,H'(x,g^{*-1}(x),\Gamma(x),p^*(x)),\,\,\forall x\in \Omega_{x_0x_{1}}.\leqno(Opt)$$
Finally, the boundary conditions
$$n_{k}p^{\ast k}_{ij}|_{\partial\Omega_{x_0x_{1}}}=\left[\frac{\partial \chi}{\partial g^{ij}}+\frac{\partial \chi}{\partial g^{ji}}\right]_{\partial\Omega_{x_0x_{1}}} \leqno (\partial\Omega_{x_0x_{1}})$$
are satisfied.}


\section{Flux-type optimal control problems }

Throughout this section, the basic geometric ingredients have the same 
significance as above; that is,  $(M,g)$  is a Riemannian manifold with local coordinates $x=(x^1,...,x^n)$ and 
$\Omega_{x_{0}x_{1}}$ denotes a subset of $M$ diffeomorphic with
a hyper-parallelepiped in $R^n$. The main goal of the section consists in 
analyzing two {\it flux-type Riemannian optimal control problems}, resulting in {\it bang-bang-type optimal solutions}.  
The key ideea is to take $x=(x^1,...,x^n)$ like an evolution (deformation) parameter.

\subsection{Optimization of total divergence}

In this subsection, $X$ is a fixed vector field on $\Omega_{x_0x_1}$. The optimal control problem we are looking to 
solve consists in finding the {\it control connection} 
$$(\Gamma^k_{ij})\in \mathcal{U}=\{\Gamma:\Omega_{x_{0}x_{1}}\rightarrow [-1, 1]^{n^3}|\,\,\Gamma_{ij}^k=\Gamma_{ji}^k\}$$
that maximize the total divergence of $X$. More precisely, we try to 
find the optimal linear connection and the optimal Riemannian structure that maximize the Bolza-type functional
$$J[\Gamma]=\int_{\Omega_{x_0x_1}}{\rm Div}\,X \,dv=\int_{\Omega_{x_0x_1}}{\rm Div}\,X \sqrt{g}\,dx,\leqno(J)$$
where $g=det(g_{ij})$, subject to the dual metric compatibility evolution PDE system
$$\frac{\pa g^{ij}}{\pa x^k}(x)=-g^{ps}(x)\left[\delta^i_p\Gamma^j_{sk}(x)+\delta^j_p\Gamma^i_{sk}(x)\right],\,\,i,j,k=1,...,n,\leqno(PDE')$$
with initial condition
$$g^{ij}(x_0)=\eta^{ij}.\leqno(x'_0)$$

{\bf Remark.} By applying the {\it Divergence Theorem} in Riemannian setting, we may rewrite the functional $J[\Gamma]$ as
$$J[\Gamma]=\int_{\pa\Omega_{x_0x_1}}g(X,N^g)\,d^g\sigma,$$
or, using local coordinates,
$$J[\Gamma]=\int_{\pa\Omega_{x_0x_1}}X^in_i\sqrt{g}\,d\sigma.$$
In the expressions above, if $N^g=(N^i)$ is the outpointing normal vector field on the boundary, with respect to the metric 
$g$, then $n=(n_i=g_{ij}N^j)$, denotes the Euclidean normal covector along the boundary.
This allows us to identify the running cost, respectively the boundary cost associated to this optimal control problem:
 $$X(x,g^{-1},\Gamma)=0;\hspace{1cm}\chi(x,g^{-1})=X^i\sqrt{g}\,n_i,$$
giving the corresponding control Hamiltonian
$$H'(x,g^{-1},\Gamma,p)=-g^{is}\Gamma^j_{sk}p^k_{ij}.$$
Since this Hamiltonian is linear with respect to the control components
$\Gamma^k_{ij}$, we have no interior optimal control $\Gamma^k_{ij}$; for optimum, 
the control must be at a vertex of $[-1,1]^{n^3}$ (see linear optimization, simplex method).

Writing the adjoint PDE system
$$\frac{\pa p^{k}_{ij}}{\pa x^k}=p^{k}_{ls}\left[\delta^l_i\Gamma^{s}_{jk}+\delta^l_j\Gamma^{s}_{ik}\right],$$
we obtain the immediate solution $p^{*k}_{ij}=C^{k}g^{*}_{ij},$ with $C=(C^k):\Omega_{x_0x_1}\rightarrow R^n$, $\frac{\pa C^k}{\pa x^k}=0.$
Then
$$H'(x,g^{*-1},\Gamma,p^*)=-C^k\Gamma^s_{ks}.$$
Therefore, the optimal control maximizing the total divergence is a linear connection having the bang-bang-type components
  $$\Gamma^{*k}_{ij}=\left\{\begin{array}{lll}
\delta_{jl} \epsilon^l & \mbox{if} & k=i,\, \epsilon^j\neq0\\
\delta_{il} \epsilon^l & \mbox{if} & k=j,\, \epsilon^i\neq0\\
 \mbox{arbitrary}, & \mbox{otherwise,}\\
 \end{array}\right.$$
 where $\epsilon ^l=\sgn(-C^l).$
 
 Moreover, the boundary constraints corresponding to these solutions are
$$[n_kC^k g^{*}_{ij}](x)=[n_k \left(\sqrt{g^{*}}X^k\right)g^{*}_{ij}](x),\,\,\forall x\in\pa\Omega_{x_0x_1},$$
that is 
$$n_k(x)\left(C^k-\sqrt{g^{*}}X^k\right)(x)=0,\,\,\forall x\in\pa\Omega_{x_0x_1}$$
and, together with the initial condition
$$g^{-1}(x_0)=\eta$$
may help us to determine the solenoidal vector field $C$.

{\bf Remark.} If we replace the above maximum-type problem with a minimizing one,
 we obtain similar solutions, with $\epsilon ^l=\sgn(C^l),\,\,l=1,...,n.$

\subsection{Optimization of total Laplacian}

Let $f:\Omega_{x_0x_1}\rightarrow R$ be a fixed differentiable function.
The optimal control problem we are interested in consists in maximizing the functional
$$J[\Gamma]=\int_{\Omega_{x_0x_1}}\Delta^gf dv=\int_{\pa\Omega_{x_0x_1}}g^{ij}f_i n_j\sqrt{g} d\sigma,\leqno(J)$$
subject to the dual metric compatibility evolution PDEs system
$$\frac{\pa g^{ij}}{\pa x^k}(x)=-g^{ps}(x)\left[\delta^i_p\Gamma^j_{sk}(x)+\delta^j_p\Gamma^i_{sk}(x)\right],\,\,i,j,k=1,...,n,\leqno(PDE')$$
with control restriction
$$(\Gamma^k_{ij})\in \mathcal{U}=\{\Gamma:\Omega_{x_{0}x_{1}}\rightarrow [-1, 1]^{n^3}|\,\,\Gamma_{ij}^k=\Gamma_{ji}^k\}$$
and with initial condition
$$g^{ij}(x_0)=\eta^{ij}.\leqno(x'_0)$$

The running cost, respectively the boundary cost associated to this optimal control problem are

$$X(x,g^{-1},\Gamma)=0;\hspace{1cm}\chi(x,g^{-1})=g^{ij}\sqrt{g}\,f_in_{j},$$
where $f_k(x)=\frac{\pa f}{\pa x^k}$ and $g^{ij}=g^{ij}(x)$ denote the components of the inverse metric matrix. Again, the control Hamiltonian is
$$H'(x,g^{-1},\Gamma,p)=-g^{is}\Gamma^j_{sk}p^k_{ij}$$
and, since it is linear with respect to the control components
$\Gamma^k_{ij}$, we have no interior optimal control $\Gamma^k_{ij}$.

Writing the Riemannian maximum principle gives us the adjoint PDE system
$$\frac{\pa p^{k}_{ij}}{\pa x^k}=p^{k}_{ls}\left[\delta^l_i\Gamma^{s}_{jk}+\delta^l_j\Gamma^{s}_{ik}\right],$$
with same possible solution as in the previous section, that is $p^{*k}_{ij}=C^{k}g^{*}_{ij},$ 
with $C=(C^k):\Omega_{x_0x_1}\rightarrow R^n$, $\frac{\pa C^k}{\pa x^k}=0.$
Replacing within the control Hamiltonian, we obtain
$$H'(x,g^{*-1},\Gamma,p^*)=-C^a\Gamma^s_{as}.$$
Therefore, the optimal control maximizing the gradient flux is of bang-bang-type 
  $$\Gamma^{*k}_{ij}=\left\{\begin{array}{lll}
\delta_{jl} \epsilon^l & \mbox{if} & k=i,\, \epsilon^j\neq0\\
\delta_{il} \epsilon^l & \mbox{if} & k=j,\, \epsilon^i\neq0\\
 \mbox{arbitrary}, & \mbox{otherwise,}\\
 \end{array}\right.$$
 where $\epsilon ^l= \sgn (-C^l).$

This time instead, the boundary constraints generating the solenoidal tensor field $C$ are
$$[n_kC^k g^{*}_{ij}](x)=\left[n_k g^{*kl}\left(f_ig^*_{lj}+f_jg^*_{li}-f_lg^*_{ij}\right)\sqrt{g^{*}}\right](x),\,\,\forall x\in\pa\Omega_{x_0x_1}.$$

{\bf Remark.} For both the foregoing problems, we may look for some particular solutions.
\begin{enumerate}
\item We may chose
$${\Gamma^*}^{k}_{ij}=\delta^k_i(\delta_{jl}\epsilon^l)+\delta^k_j(\delta_{il}\epsilon^l)-\delta_{ij}\delta^{kp}(\delta_{pl}\epsilon^l),$$
that is $\Gamma^*$ is an Euclidean conformal linear connection (see \cite{[5]}), i.e. $\Gamma^*$ is
conformal with the Levi-Civita connection associated to the Euclidean metric $g^0_{ij}=\delta_{ij}$.
Then, the optimal Riemannian metric 
$$g^{*ij}=K\delta^{ij}e^{-2\delta_{kl}\epsilon^kx^l}$$
is a {\it soliton-type solution} for the dual metric compatibility evolution $(PDE)$
$$\frac{\pa g^{ij}}{\pa x^k}(x)=-g^{ps}(x)\left[\delta^i_p\Gamma^{*j}_{sk}(x)+\delta^j_p\Gamma^{*i}_{sk}(x)\right]$$
and, also, is a dual Riemannian structure (a dual Riemannian metric having  $\Gamma^*$ as Levi-Civita connection).

\item We may consider $$\Gamma^k_{ij}=\epsilon^k\epsilon_i\epsilon_j,$$
where $\epsilon_i=\delta_{ij}\epsilon^j$.
Then, the $(PDE')$ system writes
$$\frac{\pa g^{ij}}{\pa x^k}=\left(g^{is}\epsilon^j+g^{js}\epsilon^i\right)\epsilon_s\epsilon_k,$$
admitting the following soliton-type solution
$$g^{ij}=\left[\alpha e^{-2n\epsilon_kx^k}+\frac{\alpha^i+\alpha^j}{2}e^{-n\epsilon_kx^k}\right]\epsilon^i\epsilon^j \,\,\,\mbox{(no summation)},$$
where $\alpha,\,\,\alpha^i$ denote real constants, satisfying $\sum_{i=1}^n \alpha^i=0.$
The disadvantage of the latter solution is that it may not be a Riemannian 
structure, but only a symmetric $(2,0)$-type tensor field or, in best case scenario 
(i.e. $\epsilon^i\neq0,\,\,\forall i=1,...,n$), a semi-Riemannian structure. 
\end{enumerate}

\section{The optimal geometry of pipes}

This section is meant to emphasize the practical utility of the theoretical facts 
described above, by analyzing a classical problem in Hydraulics and Fluid Mechanics. 
Given a pipe, in the general sense (water pipe, gas pipe, blood vessel) containing a fluid flow, 
it is well known that the Divergence Theorem helps us to measure the flux of the fluid flow 
through pipe walls. Sometimes instead, for practical reasons, it is of major utility to identify 
the optimal shape of the pipe, allowing the minimum flux through walls. This is the problem 
analyzed in this section. More precisely, given the directionality of the fluid through the pipe, 
we decide about the best way to conceive the pipe (the optimal geometric shape), 
such that the flux of the fluid through pipe walls to be minimal.

Let $D^1$ denote the closed disc of radius one and let $M=D^1\times (0,1)$ be a differential 
manifold with boundary describing the interior and the boundary of a cylinder. We identify the pipe in the 
Euclidean space (in the sense of some diffeomorphism) with the manifold $M$. Given a 
flow through the pipe, described by a vector field $F=X\frac{\pa}{\pa x}+Y\frac{\pa}{\pa y}+Z\frac{\pa}{\pa z}$ 
on $M$, we shall find a Riemannian structure on $M$, minimizing the flux of $F$.  For this, 
we consider the local map $V=M-\{(x,y,z)\in M|\,\,y=0,\,\,x\geq0\}$. Using the cylindrical coordinates 
$(\rho,\theta,z)$, we may identify $V$ with the parallelepiped $(0,1]\times(0,2\pi)\times(0,1)$ in $R^3$. 
Moreover, we suppose that the expression of $F$ with respect to these new 
coordinates is $F=R\frac{\pa}{\pa \rho}+T\frac{\pa}{\pa\theta}+\zeta\frac{\pa}{\pa z}$. Then,
$$\left\{\begin{array}{lll}
R(\rho,\theta,z)=X(\rho \cos\theta, \rho \sin\theta,z)\cos\theta+Y(\rho \cos\theta, \rho \sin\theta,z)\sin\theta;\\
T(\rho,\theta,z)=\rho\left[-X(\rho \cos\theta, \rho \sin\theta,z)\sin\theta+Y(\rho \cos\theta, \rho \sin\theta,z)\cos\theta\right];\\
\zeta(\rho,\theta,z)=Z(\rho \cos\theta, \rho \sin\theta,z),
\end{array}\right.$$
or, conversely,
$$\left\{\begin{array}{lll}
X(x,y,z)=R(\sqrt{x^2+y^2},\arctan{\frac{y}{x}},z)\frac{x}{\sqrt{x^2+y^2}}-T(\sqrt{x^2+y^2},\arctan{\frac{y}{x}},z)y;\\
Y(x,y,z)=R(\sqrt{x^2+y^2},\arctan{\frac{y}{x}},z)\frac{y}{\sqrt{x^2+y^2}}+T(\sqrt{x^2+y^2},\arctan{\frac{y}{x}},z)x;\\
Z(x,y,z)=\zeta(\sqrt{x^2+y^2},\arctan{\frac{y}{x}},z).
\end{array}\right.$$

Applying the results derived in the previous section, we obtain the optimal  Euclidean conformal structure
$$g(\rho,\theta,z)=Ke^{2sgn(R(1,\theta,z))\rho}\left(d\rho^2+d\theta^2+dz^2\right).$$
Using the above relations between the components of $F$ relative to the cylindrical and Cartesian 
coordinates we derive also the Cartesian expression of the optimal Riemannian structure: 
$$g=Ke^{2S\left(\frac{x}{\sqrt{x^2+y^2}},\frac{x}{\sqrt{x^2+y^2}},z\right)\sqrt{x^2+y^2}}\left(\begin{array}{cccl}1 &0 &0\\
                                             0 &x^2+y^2 &0\\
                                             0& 0& 1\end{array}\right),$$
where $S=sgn\langle N, F\rangle$ on the boundary,
$\langle\,, \rangle$ denotes the canonical inner product  on $R^3$ and $N(x,y,z)=x\frac{\pa}{\pa x}+y\frac{\pa}{\pa y}$ 
is the normal vector field along the boundary of $M$. 

In conclusion, the direction of the flow $F$ as vector field in $R^3$ is directly involved in 
the geometric configuration of the pipe. More precisely, if $F$ is pointed outward,  then the 
diameter of the pipe increases; conversely, if $F$ is pointed inward, the diameter decreases. 
Therefore, the optimal shape for the pipe walls is the one tangent, at each point, to the flow $F$.

\end{document}